\documentclass[11pt]{article}
\usepackage{amsfonts}
\usepackage{graphicx,amssymb,amsmath}

\newtheorem{definition}{Definition}
\newtheorem{remark}{Remark}
\newtheorem{theorem}{Theorem}

\newtheorem{proposition}{Proposition}

\begin{document}


\begin{center}
LINEAR PROCESSES FOR FUNCTIONAL DATA

Andr\'e \textsc{Mas}$^{a}$ and Besnik \textsc{Pumo}$^{b}$\\[0pt]

$^{a}$ I3M, Universit\'e Montpellier 2,\\[0pt]
Place Eug\`ene Bataillon, 34095 Montpellier, France\\[0pt]
mas@math.univ-montp2.fr\\[0pt]
~\\[0pt]
$^{b}$ Agrocampus Ouest, Centre d'Angers,\\[0pt]
2 rue Le N\^otre, 49045 Angers, France \\[0pt]
besnik.pumo@agrocampus-ouest.fr\\[0pt]
~\\[0pt]
\end{center}

\section*{Abstract}

Linear processes on functional spaces were born about fifteen years ago. And
this original topic went through the same fast development as the other
areas of functional data modeling such as PCA or regression. They aim at
generalizing to random curves the classical ARMA models widely known in time
series analysis. They offer a wide spectrum of models suited to the
statistical inference on continuous time stochastic processes within the
paradigm of functional data. Essentially designed to improve the quality and
the range of prediction, they give birth to challenging theoretical and
applied problems. We propose here a state of the art which emphasizes recent
advances and we present some promising perspectives based on our experience
in this area.

\section{Introduction}

The aim of this chapter is double. First of all we want to provide the
reader with basic theory and application of linear processes for functional
data. The second goal consists for us in giving a state of the art which
complements the monograph by Bosq (2000). Many crucial theorems were given
in this latter book to which we will frequently refer. Consequently, even if
our work is self-contained we pay special attention to recent results,
published from 2000 to 2008, and try to draw the lines of future and
promising research in this area.

It is worth recalling now the approach that leads to modelizing and
inferring from curves-data. We start from a continuous time stochastic
process $\left( \xi_{t}\right) _{t\geq0}$. The paths of $\xi$ are cut into
equally spaced pieces of trajectories. Each of these piece is then viewed as
a random curve. With mathematical symbols we set:%
\begin{equation*}
X_{k}\left( t\right) =\xi_{kT+t},\quad0\leq t\leq T
\end{equation*}
where $T$ is fixed. The function $X_{k}\left( \cdot\right) $ maps $\left[ 0,T%
\right] $ to $\mathbb{R}$ and is random. Observing $\xi$ over $\left[ 0,nT%
\right] $ produces a $n$ sample $X_{1},...,X_{n}$. Obviously the choice of $%
T $ is crucial and is usually left to the practitioner and may be linked
with seasonality (with period $T$). Dependence along the paths of $\xi$ will
create dependence between the $X_{i}$'s. But this approach is not restricted
to the whole path of stochastic process. One could as well imagine to model
whole curves observed at discrete intervals: the interest rate curves at day 
$k$ $I_{k}\left( \delta\right) $ is for instance a function linking duration 
$\delta$ (as an input) and the associated interest rates (as outputs).
Observing these curves, whose random variations will depend on financial
markets, along $n$ days produces a sample similar in nature to the one
described above, although there is no underlying continuous-time process in
this situation, rather a surface $\left(k, \delta, I_{k}(\delta) \right)$.
We refer for instance to Kargin and Onatski (2008) for an illustration.

Statistical models will then be proposed and mimic or adapt the scalar or
finite-dimensional approaches for time-series (see Brockwell, Davis (1987)).
Each of these (random or not) functions will be viewed as a vector in a
vector space of functions. This paradigm has been adopted for a long time in
probability theory as will be seen through, for instance, Ledoux and
Talagrand (1991) and references therein. But the first book entirely
dedicated to the formal and applied aspects of statistical inference in this
setting is certainly due to Ramsay and Silverman (1997), followed by Bosq
(2000), Ramsay and Silverman (2002) again then Ferraty and Vieu (2006).

In the sequel we will consider centered processes with values in a Hilbert
space of functions denoted $H$ with inner product $\left\langle \cdot ,\cdot
\right\rangle $ and norm $\left\Vert \cdot \right\Vert $. The Banach setting
though more general has several drawbacks. Some references will be given yet
throughout this section. The reason for privileging Hilbert spaces are both
theoretic and practical. First many fundamental asymptotic theorems are
stated under simple assumptions in this setting. The central limit theorem
is a good example. Considering random variables with values in $C\left( %
\left[ 0,1\right] \right) $ or in H\"{o}lder spaces for instance lead to
very specific assumptions to get the CLT and computations are often uneasy
whereas in a Hilbert space moment conditions are usually both necessary and
sufficient. The nice geometric features of Hilbert space allow us to
consider denumerable bases, projections, etc in a framework that generalizes
the euclidean space with few drawbacks. Besides, in practice, recovering
curves from discretized observations by interpolation or smoothing
techniques such as splines or wavalets yields functions in the Sobolev
spaces, say $W^{m,2}$ (here $m$ is an order of differentiation connected
with the desired smoothness of the output), are all Hilbert spaces. We refer
to Ziemer (1989) or to Adams, Fournier (2003) for monographs on Sobolev
spaces.

In statistical models, unknown parameters will be functions or linear
operators (the counterpart of matrices of the euclidean space), the latter
being of utter interest. We give now some basic facts about operators which
will be of great use in the sequel.

Several monographs are dedicated to operator theory, which is a major theme
within the mathematical science. Classical references are Dunford, Schwartz
(1988) and Gohberg, Goldberg and Kaashoek (1991){\scriptsize . }The adjoint
of the operator $T$ is classically denoted $T^{\ast}$. The Banach space of
compact operators $\mathcal{C}$ on a Hilbert space $H$ is separable when
endowed with the classical operator norm $\left\Vert \cdot\right\Vert
_{\infty}$:%
\begin{equation*}
\left\Vert T\right\Vert _{\infty}=\sup_{x\in\mathcal{B}_{1}}\left\Vert
Tx\right\Vert
\end{equation*}
where $\mathcal{B}_{1}$ denotes the unit ball of the Hilbert space $H$. The
space $\mathcal{C}$ contains the set of Hilbert-Schmidt operators which is a
Hilbert space and denoted $\mathcal{S}$. Let $T$ and $S$ belong to $\mathcal{%
S}$ the inner product between $T$ and $S$ and the norm of $T$ are
respectively defined by:%
\begin{align*}
\left\langle S,T\right\rangle _{\mathcal{S}} & =\sum_{p}\left\langle
Se_{p},Te_{p}\right\rangle , \\
\left\Vert T\right\Vert _{\mathcal{S}}^{2} & =\sum_{p}\left\Vert
Te_{p}\right\Vert^{2}
\end{align*}
where $\left( e_{p}\right) _{p\in\mathbb{N}}$ is a complete orthonormal
system in $H.$ The inner product and the norm defined just above do not
depend on the choice of the c.o.n.s.$\left( e_{p}\right) _{p\in\mathbb{N}}$.
The nuclear (or trace-class) operators are another important family of
operators for which the series:%
\begin{equation*}
\sum_{p}\left\Vert Te_{p}\right\Vert <+\infty.
\end{equation*}
It is plain that a trace class operator is Hilbert-Schmidt as well. Many of
the asymptotic result mentioned from now on and involving random operators
are usually obtained for the Hilbert-Schmidt norm, unless explicitly
mentioned. It should be noted as well that this norm is thinner than the
usual operator norm.

The next section is devoted to general linear processes. Then we will focus
on the autoregressive model and its recent advances, which will be developed
in the third section. We will conclude with some issues for future work.

\section{General linear processes}

The linear processes on function spaces generalize the classical scalar or
vector linear processes to random elements which are curves or functions and
more generally valued in an infinite-dimensional separable Hilbert space $H$.

\begin{definition}
Let $\left( \epsilon_{k}\right) _{k\in\mathbb{N}}$ be a sequence of i.i.d.
centered random elements in $H$ and let $\left(a_{k}\right)_{k \in \mathbb{N}%
}$ be a sequence of bounded linear operators from $H$ to $H$ such that $%
a_{0}=I$ and $\mu\in H$ be a fixed vector. If 
\begin{equation}
X_{n}=\mu + \sum_{j=0}^{+\infty}a_{j}\left( \epsilon_{n-j}\right) ,
\label{AMBP_def-lin-proc}
\end{equation}
$\left( X_{n}\right) _{n\in\mathbb{N}}$ is a linear process on $H$ (denoted
in the sequel $H$-linear process) with mean $\mu$.
\end{definition}

Unless explicitly mentioned the mean function $\mu$ will always be assumed
to be null (and the process $X$ is centered). Its seems that, after a
collection of paper dating back to the end of the 90's-early 00's the model
creates less inspiration in the community. We guess that the recent works by
Bosq (2007) and the book by Bosq and Blanke (2007) may bring some fresh
ideas. We state here some basic facts: invertibility and convergence of
estimated moments.

\subsection{Invertibility}

When the sequence $\epsilon $ is a strong $H$-white noise, that is a
sequence of i.i.d. random elements such that $E\left\Vert \epsilon
\right\Vert ^{2}<+\infty $ and whenever 
\begin{equation}
\sum_{j=0}^{+\infty }\left\Vert a_{j}\right\Vert _{\infty }^{2}<+\infty 
\label{stathyp}
\end{equation}
the series defining the process $\left( X_{n}\right) _{n\in \mathbb{N}}$
through (\ref{AMBP_def-lin-proc}) converges in square norm and almost surely
through the $0-1$ law. The strict stationarity of $X_{n}$ is ensured as
well. The problem of invertibility is addressed in Merlev\`{e}de (1995).

\begin{theorem}
If $\left( X_{n}\right) _{n\in \mathbb{N}}$ is a linear process with values
in $H$ defined by (\ref{AMBP_def-lin-proc}) and such that :%
\begin{equation}
1-\sum_{j=1}^{+\infty }z^{j}\left\Vert a_{j}\right\Vert _{\infty }\neq
0\quad for\mathrm{\ }\left\vert z\right\vert <1  \label{AMBP_cond_inv_lph}
\end{equation}%
then $\left( X_{n}\right) _{n\in \mathbb{N}}$ is invertible:%
\begin{equation*}
X_{n}=\epsilon _{n}+\sum_{j=1}^{+\infty }\rho _{j}\left( X_{n-j}\right) 
\end{equation*}%
where all the $\rho _{j}$'s are bounded linear operators in $H$ with $%
\sum_{j=1}^{+\infty }\left\Vert \rho _{j}\right\Vert _{\infty }<+\infty $
and the series converges in mean square and almost surely.
\end{theorem}

\begin{remark}
We deduce from the latter that $\epsilon_{n}$ is the innovation of the
process $X$ and that (\ref{AMBP_def-lin-proc}) coincides with the Wold
decomposition of $X.$
\end{remark}

We give now some convergence theorems for the mean and the covariance of
Hilbert-valued linear processes. These results are not completely new but
essential.

\subsection{Asymptotics}

It is worth mentioning a general scheme for proving asymptotic results for
linear processes. If several approaches are possible, it turns out that, up
to the authors' opinion, one of the most fruitful relies on approximating
the process $X_{n}$ by truncated versions like:%
\begin{equation*}
X_{n,m}=\sum_{j=0}^{m}a_{j}\left( \epsilon_{n-j}\right)
\end{equation*}
where $m\in\mathbb{N}$. The sequence $X_{n,m}$ is for fixed $m,$ blockwise
independent: $X_{n+m+1,m}$ is indeed stochastically independent from $%
X_{n,m} $ if the $\epsilon_{j}$'s are. The outline of the proofs usually
consists in proving asymptotic results for the $m$-dependent sequence $%
X_{n,m},$ then to let $m$ tend to infinity with an accurate control of the
residual $X_{n}-X_{n,m}=\sum_{j=m+1}^{+\infty}a_{j}\left( \epsilon
_{n-j}\right) .$

\subsubsection{Mean}

Asymptotic results for the mean of a linear process may be found in Merlev%
\`{e}de (1996) and Merlev\`{e}de, Peligrad and Utev (1997). Even if the
first is in a way more general, we focus here on the second article since it
deals directly with the mean of the non-causal process indexed by $\mathbb{Z}
$ :%
\begin{equation*}
X_{k}=\sum_{j=-\infty}^{+\infty}a_{j}\left( \epsilon_{k-j}\right) .
\end{equation*}
The authors obtain sharp conditions for the CLT of $S_{n}=%
\sum_{k=1}^{n}X_{k}.$

\begin{theorem}
Let $\left( a_{j}\right) _{j\in\mathbb{Z}}$ be a sequence of operators such
that:%
\begin{equation*}
\sum_{j=-\infty}^{+\infty}\left\Vert a_{j}\right\Vert _{\infty}^{2}<+\infty
\end{equation*}
Then%
\begin{equation*}
\frac{S_{n}}{\sqrt{n}}\rightarrow_{w}N\left( 0,A\Gamma_{\epsilon}A^{\ast
}\right)
\end{equation*}
where $N\left( 0,A\Gamma_{\epsilon}A^{\ast}\right) $ is the $H$-valued
centered gaussian random element with covariance operator $A\Gamma
_{\epsilon}A^{\ast}$ where $\Gamma_{\epsilon}=E\left( \epsilon
_{0}\otimes\epsilon_{0}\right) $ is the covariance operator of $\epsilon_{0}$
and $A=\sum_{j=-\infty}^{+\infty}a_{j}$.
\end{theorem}

Remind that if $u$ and $v$ are two vectors in $H$ then notation $u\otimes v$
stands for the rank-one linear operator from $H$ to $H$ defined by: $\left(
u\otimes v\right) \left( x\right) =\left\langle v,x\right\rangle u$.

This result is extended with additional assumptions to the case of strongly
mixing $\epsilon _{k}$'s. Note that the problem of weak convergence for the
mean of stationary Hilbertian process under mixing conditions had been
addressed in the early 80's by Maltsev and Ostrovski (1982). A standard
equi-integrability argument and classical techniques provide the following
rates of convergence for $S_{n}$. Nazarova (2000) proved the same sort of
theorem when $X$ is a linear random field with values in a Hilbert space.
Now we turn to the rate of convergence of the empirical mean in quadratic
mean and almost surely. The following theorem may be found in Bosq (2000).

\begin{proposition}
Let $X_{k}=\sum_{k=0}^{+\infty}a_{j}\left( \epsilon_{k-j}\right) $ and $%
S_{n}=\sum_{k=1}^{n}X_{k}$ then 
\begin{align*}
nE\left\Vert \frac{S_{n}}{n}\right\Vert ^{2} & \rightarrow\sum_{k=-\infty
}^{+\infty}E\left\langle X_{0},X_{k}\right\rangle , \\
\frac{n^{1/4}}{\left( \log n\right) ^{1/2+\epsilon}}\left\Vert \frac{S_{n}}{n%
}\right\Vert & \rightarrow0\quad a.s.
\end{align*}
for all $\epsilon>0.$
\end{proposition}

We turn to covariance operators now.

\subsubsection{Covariance operators}

The situation is slightly more complicated than for the mean due to the
tensor product.

\begin{definition}
The theoretical covariance operator at lag $h\in \mathbb{N}$ of a process $X$
is defined by:%
\begin{equation*}
\Gamma _{h}=\mathbb{E}\left( X_{h}\otimes X_{0}\right) .
\end{equation*}%
The linear operator $\Gamma _{h}$ is nuclear on $H$ when the second order
strong moments of the $X$ is convergent. Its empirical counterpart based on
the sample is:%
\begin{equation*}
\Gamma _{n,h}=\dfrac{1}{n}\sum_{t=1}^{n}X_{t+h}\otimes X_{t}.
\end{equation*}%
The covariance operator of the process, $\Gamma _{0}=\Gamma $ is
selfadjoint, positive and nuclear hence Hilbert-Schmidt and compact.
\end{definition}

It should be noted that $\Gamma _{h}$ is not in general a symmetric operator
conversely to the classical covariance operator $\Gamma _{0}$. The weak
convergence of covariance operators for $H$-linear processes was addressed
by Mas (2002). It is assumed that:%
\begin{align*}
\mathbb{E}\left\Vert \epsilon _{0}\right\Vert ^{4}& <+\infty \\
\sum_{k=-\infty }^{+\infty }\left\Vert a_{k}\right\Vert _{\infty }& <+\infty
\end{align*}%
then the vector of the $h$ covariance operators up to any fixed lag $h$ is
asymptotically gaussian in the Hilbert-Schmidt norm.

\begin{theorem}
\label{AMBP_main}Let us consider the following linear and Hilbert space
valued process 
\begin{equation*}
X_{t}=\sum_{j=-\infty}^{+\infty}a_{j}\left( \epsilon_{t-j}\right)
\end{equation*}
then 
\begin{equation*}
\sqrt{n}\left( 
\begin{array}{c}
\Gamma_{n,0}-\Gamma_{0} \\ 
\Gamma_{n,1}-\Gamma_{1} \\ 
... \\ 
\Gamma_{n,h}-\Gamma_{h}%
\end{array}
\right) \underset{n\rightarrow+\infty}{\overset{w}{\rightarrow}}G_{\Gamma}
\end{equation*}
where $G_{\Gamma}=\left( G_{\Gamma}^{(0)},...,G_{\Gamma}^{(h)}\right) $ is a
Gaussian centered random element with values in $\mathcal{S}^{h+1}.$ Its
covariance operator is $\Theta_{\Gamma}=\left(
\Theta_{\Gamma}^{(p,q)}\right) _{0\leq p,q\leq h}$ which is a nuclear
operator in $\mathcal{S}^{h+1}$ defined blockwise for all $T$ in $\mathcal{S}
$ by 
\begin{equation}
\Theta_{\Gamma}^{(p,q)}\left( T\right) =\sum_{h}\Gamma_{h+p-q}T\Gamma
_{h}+\sum_{h}\Gamma_{h+q}T\Gamma_{h-p}+A_{q}\left( \Lambda-\Phi\right)
A_{p}\left( T\right)  \label{AMBP_lafin}
\end{equation}
where $\Lambda,$ $\Phi$ and $A_{p}$ are linear operators from $\mathcal{S}$
to $\mathcal{S}$ respectively defined by 
\begin{align*}
\Lambda\left( T\right) & =\mathbb{E}\left( \left(
\epsilon_{0}\otimes\epsilon_{0}\right) \widetilde{\otimes}\left( \epsilon
_{0}\otimes\epsilon_{0}\right) \right) \left( T\right) \\
\Phi\left( T\right) & =C\left( T+T^{\ast}\right) C+\left( C\widetilde{\otimes%
}C\right) \left( T\right)
\end{align*}
and $A_{p}\left( T\right) =\sum_{i}a_{i+p}Ta_{i}^{\ast}.$
\end{theorem}

As by-products, weak convergence results for the eigenelements of $%
\Gamma_{n,0}-\Gamma_{0}$, that is for the PCA of the stationary process $X$,
are derived. The reader interested with these developments should refer to
the paper by Mas and Menneteau (2003a) which proposes a general method to
derive asymptotics for the eigenvalues and eigenvectors of $\Gamma_{n,0}$
(the by-products of the functional PCA) from the covariance sequence itself.
Perturbation theory is the main tool through a modified delta-method.

\subsection{Perspectives and trends: towards generalized linear processes ?}

\label{AMBP_LHP_perspectives}

It turns out that the literature based on inference methods for general
linear processes is rather meager. Obviously estimating simultaneously many $%
a_{j}$'s seems to be intricate and not necessarily needed as the functional
AR process, which will be enlightened in the next section, is quite
successful and easier to handle. However general linear processes are the
starting point for very interesting theoretical problems where dependence
plays a key-role. We mention at last the abstract papers by Merlev\`{e}de
and Dedecker (2003) especially section 2.4 dedicated to proving a
conditional central limit theorem for linear processes under mild
assumptions and to Merlev\`{e}de and Dedecker (2007) whose section 3 deals
with rates in the law of large numbers. These works may provide theoretical
material to go further into the asymptotic study of these processes.

In a very recent article, Bosq (2007) introduces the notion of linear
process in the wide sense. The definition remains essentially the same as in
display (\ref{AMBP_def-lin-proc}) but the operators $\left( a_{j}\right)
_{j\in \mathbf{N}}$ may then be unbounded; which finally generalizes the
notion. A key role is played by linear closed spaces (LCS) introduced by
Fortet. A LCS $\mathcal{G}$ is a subspace of $L_{H}^{2}$ -the space of
random variables with values in $H$ and finite strong second moment- such
that:

\begin{description}
\item[(i)] $\mathcal{G}$ is closed in $H.$

\item[(ii)] If $X\in\mathcal{G}$, $l\left( X\right) \in\mathcal{G}$ for all
bounded linear operator $l$.
\end{description}

This theory -involving projection on LCS, weak and strong orthogonality,
dominance of operators- allows Bosq to revisit and extend the notions of
linear process, Wold decomposition, Markovian process when the bounded $%
\left( a_{j}\right) _{j\in \mathbf{N}}$ may be replaced with measurable
mappings $\left( l_{j}\right) _{j\in \mathbf{N}}$. Several examples are
given : derivatives of functional processes like in the MAH $X_{n}=\epsilon
_{n}+c\epsilon _{n}^{\prime }$, arrays of linear processes, truncated
Ornstein-Uhlenbeck process... The personnal communication Bosq (2009)
discusses these extensions to tensor products of linear processes and will
certainly shed a new light at their covariance structure. We also refer to
chapters 10 and 11 in the book by Bosq and Blanke (2007) for an exposition
of these concepts.

\section{Autoregressive processes}

\subsection{Introduction}

The model generalizes the classical AR(1) for scalar or multivariate time
series to functional data and was introduced for the first time in Bosq
(1991). Let $X_{1},\ldots ,X_{n}$ be a sample of random curves for which a
stochastic dependence is suspected (for instance the curve of temperature
observed during $n$ days at a given place). We assume that all the $X_{i}$'s
are valued in a Hilbert space $H$ and set:%
\begin{equation}
X_{n}=\rho \left( X_{n-1}\right) +\epsilon _{n}  \label{AMBP_ARH}
\end{equation}%
where $\rho $ is a linear operator from $H$ to $H$ and $\left( \epsilon
_{n}\right) _{n\in \mathbb{N}}$ is a sequence of $H$ valued centered random
elements usually with common covariance operator. The model is simple, with
a single unknown operator, leaving however the possibility to decline
various assumptions either on the operator $\rho $ (linear, compact,
Hilbert-Schmidt, symmetric or not, etc) or on the dependence between the $%
\epsilon _{n}$'s. The latter are quite often independent and identically
distributed but alternatives are possible (mixing or more naturally
martingale differences). Bosq (2000) proved that assumption (\ref{stathyp})
comes down actually to the existence of $a>0$, $b\in \left[ 0,1\right[ $such
that for all $p\in \mathbb{N}$ :%
\begin{equation*}
\left\Vert \rho ^{p}\right\Vert _{\infty }\leq ab^{p}
\end{equation*}%
which ensures that (\ref{AMBP_ARH}) admits a unique stationary solution. The
process $\left( X_{n}\right) _{n\in \mathbb{N}}$ is Markov as soon as $%
E\left( \epsilon _{n}|X_{n-1},\ldots ,X_{1}\right) =0.$ As often noted the
interest of the model relies in its predictive power. The estimation of $%
\rho $ is usually the first and necessary step before deriving the
statistical predictor given the new input $X_{n+1}$: $\widehat{\rho }\left(
X_{n}\right) $. The prediction are often compared with ARMA model or with
non-parametric smoothing techniques. The global treatment of the trajectory
as a function often ensures better long-run prediction but at the expense of
more tedious numerical procedures.

\subsubsection{Representation of stochastic processes by functional AR}

\begin{figure}[h]
\centering
\includegraphics[scale=0.44]{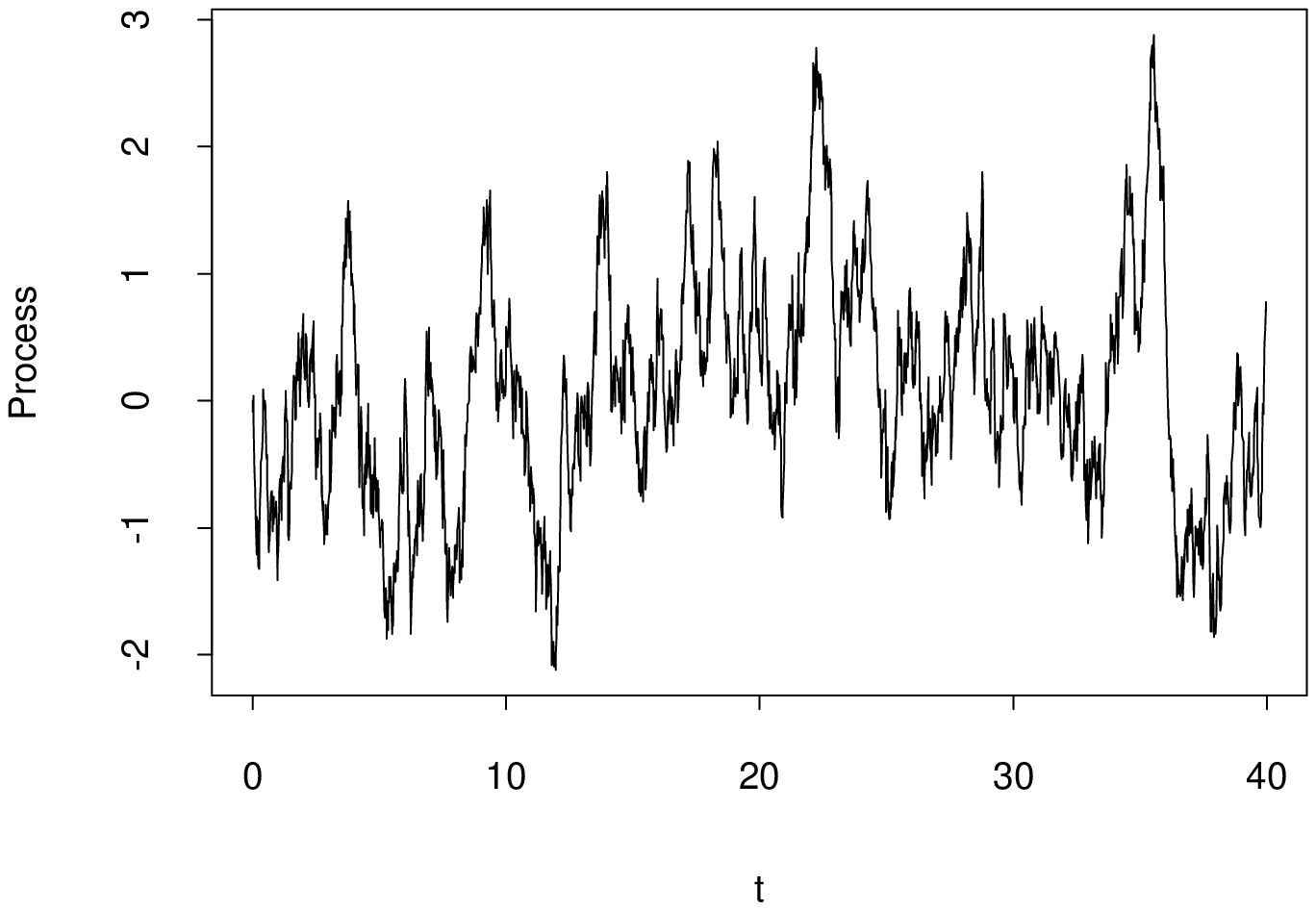} %
\includegraphics[scale=0.44]{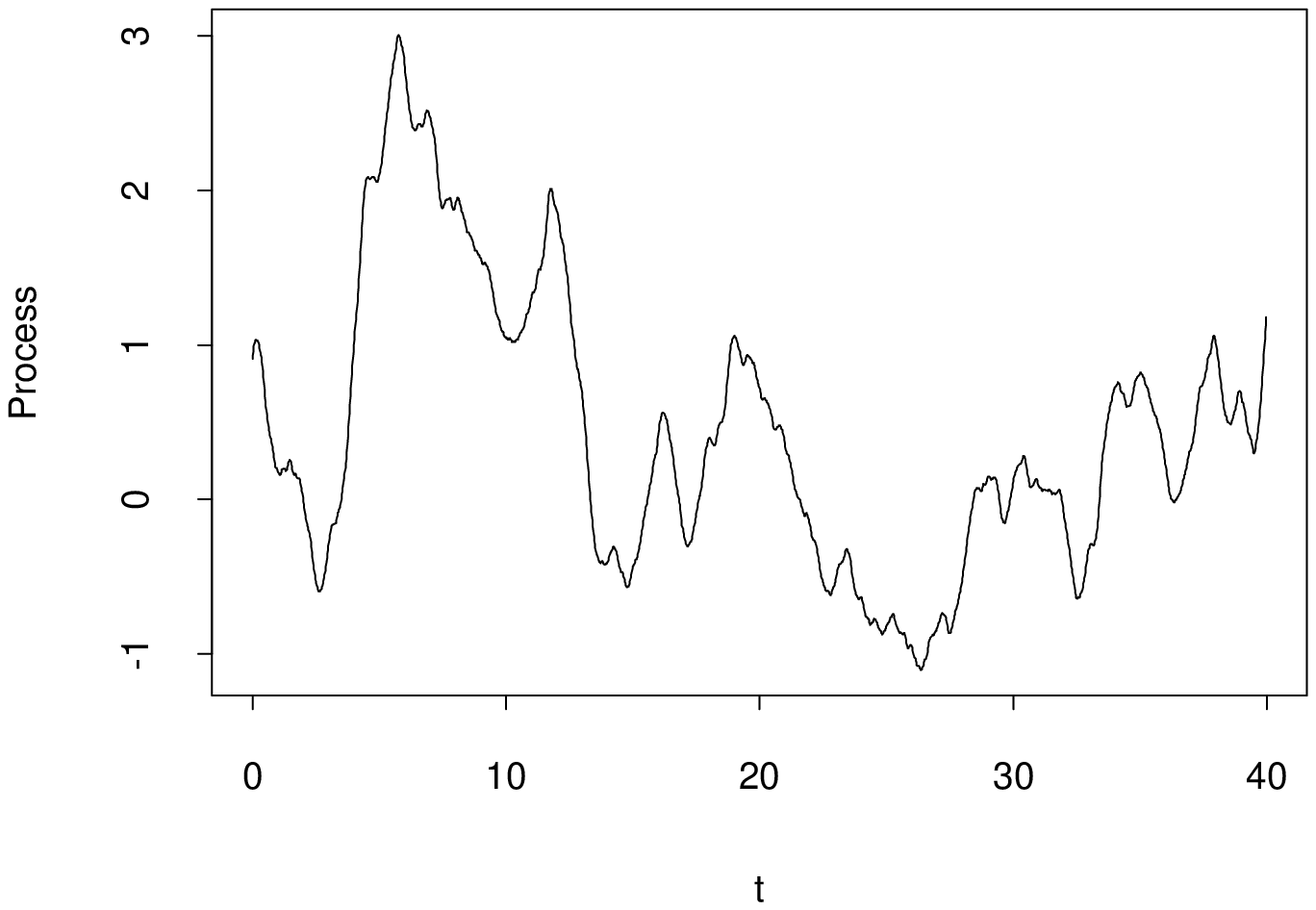}
\caption{Ornstein-Uhlenbeck ($a=1$) and Wong processes both evaluated at
instants $t_i=0.02*i$.}
\label{AMBPouprc.wgprc}
\end{figure}

Various real valued processes allow the ARH representation. We plotted on
Figure \ref{AMBPouprc.wgprc} graphs of two simulated processes, the
Ornstein-Uhlenbeck (O-U) process and the Wong process. The O-U process $%
(\eta _{t},t\in R)$ is a real stationary Gaussian process :%
\begin{equation*}
\eta _{t}=\int_{-\infty }^{t}e^{-a(t-u)}dw_{u},t\in R
\end{equation*}%
where $(w_{t})_{t\in R}$ is a bilateral standard Wiener process and $a$ a
positive constant. Bosq (1996) gives the ARH representation $X_{n}=\rho
(X_{n-1})+\epsilon _{n}$ with values in $L^{2}:=L^{2}[0,1]$ where $%
X_{n}(t)=\eta _{n+t},t\in \lbrack 0,1],n\in Z$ and $\rho $ is a degenerated
linear operator 
\begin{equation*}
\rho (x)(t)=e^{-at}x(1),t\in \lbrack 0,1],x\in L^{2}
\end{equation*}%
and 
\begin{equation*}
\epsilon _{n}(t)=\int_{n}^{n+t}e^{-a(n+t-s)}dw_{s},t\in \lbrack 0,1],n\in Z.
\end{equation*}%
The Wong process is a mean-square differentiable stationary Gaussian process
which is zero-mean and is defined for $t\in R$ by: 
\begin{equation*}
\xi _{t}=\sqrt{3}\,\exp \left( -\sqrt{3}\,t\right) \int_{0}^{\exp \left( 2t/%
\sqrt{3}\right) }w_{u}\,du.
\end{equation*}%
Cutting $R$ in intervals of length 1 and defining $X_{n}(t)=\xi _{n+t}$ for $%
t\in \lbrack 0,1]$, Mas and Pumo (2007) obtain ARH representation, $%
X_{n}=A(X_{n-1})+\epsilon _{n}$ of this process with values in Sobolev space 
$W:=W^{2,1}=\left\{ u\in L^{2},u^{\prime }\in L^{2}\right\} $ and 
\begin{equation*}
\epsilon _{n}(t)=\sqrt{3}\,\exp [-\sqrt{3}\,(n-1+t)]\int\limits_{\exp
[2(n-1)/\sqrt{3}]}^{\exp [2(n-1+t)/\sqrt{3}]}\left[ w_{u}-w_{\exp [2(n-1)/%
\sqrt{3}]}\right] du
\end{equation*}%
for $t\in \lbrack 0,1]$. The linear and degenerated operator A is given by $%
\phi +\Psi (D)$ where $D$ is the ordinary differential operator and 
\begin{equation*}
\lbrack \phi (f)](t)=[\exp (-\sqrt{3}t)+\sqrt{3}c(t)]f(1),\hspace{0.11in}%
[\Psi (D)(f)](t)=c(t)f^{\prime }(1).
\end{equation*}%
with $c(t)=\frac{\sqrt{3}}{2}\cdot \exp (-\sqrt{3}t)\cdot \{\exp (2t/\sqrt{3}%
)-1\}$.

Other examples are given in the paper of Bosq (1996) or in the classical
book of Bosq (2000).

A major issue would be to infer deeper links between autoregressive
functional processes and diffusion processes or stochastic differential
equations with more general autocorrelation operators. But this remains an
open question and the work by Ramsay (2000) about this topic should
certainly deserve more attention to be extended.

\subsubsection{Asymptotics for the mean and covariance}

Obviously all the results obtained for general linear processes hold for the
ARH(1): namely for the mean and the covariance operator. Some new results
are stated below -they are new essentially with respect to Bosq (2000)- and
are related with moderate deviations (Mas and Menneteau (2003b)) and laws of
the iterated logarithm (Menneteau (2003)). Let $\eta$ be a square integrable
real valued random variable. We need to introduce the following notations ($%
I_{X}$ and $J_{X}$ are functions from $H$ to $H$, $J_{\Gamma}$ is a function
from S to S and $K_{\Gamma}$ is a subset of $S$):%
\begin{align}
u_{1} & =\rho\left( X_{0}\right) \otimes\epsilon_{1}+\epsilon
_{1}\otimes\rho\left( X_{0}\right) +\epsilon_{1}\otimes\epsilon
_{1}-\Gamma_{\epsilon},  \label{AMBP_ui} \\
I_{X}\left( x\right) & =\sup_{h\in H}\left\{ \left\langle h,x-\rho\left(
x\right) \right\rangle -\mathbb{E}\exp\left\langle h,\epsilon
_{1}\right\rangle \right\} ,  \notag \\
J_{X}\left( x\right) & =\frac{1}{2}\inf\left\{ E\eta^{2}:x=E\left[
\eta\left( I_{H}-\rho\right) ^{-1}\left( \epsilon_{1}\right) \right]
\right\} ,  \notag \\
J_{\Gamma}\left( s\right) & =\frac{1}{2}\inf\left\{ E\eta^{2}:s=E\left[
\eta\left( I_{S}-R\right) ^{-1}\left( u_{1}\right) \right] \right\} ,  \notag
\\
K_{\Gamma} & =\left\{ E\left[ \eta\left( I_{S}-R\right) ^{-1}\left(
u_{1}\right) \right] :\eta\in L^{2}\left( P\right) ,E\eta^{2}\leq1\right\} 
\notag
\end{align}
and $R$ is a linear operator from $S$ to $S$ defined by $R\left( s\right)
=\rho s\rho^{\ast}$.

We refer to Dembo and Zeitouni (1993) for an exposition on large and
moderate deviations.

\begin{theorem}
The empirical mean of the ARH(1) process follows the large deviation
principle in $H$ with speed $n^{-1}$ and rate function $I_{X}$ and the
moderate deviation principle with rate function $J_{X}$.\newline
The covariance sequence of the ARH(1), $\Gamma_{n}-\Gamma$, follows the
moderate deviation principle in the space of Hilbert-Schmidt operators with
rate function $J_{\Gamma}$ and the law of the iterated logarithm with limit
set $K_{\Gamma}$.\newline
\end{theorem}

The first results obtained on the covariance sequence are in Bosq (1991) but
we mention here the interesting decomposition given in Bosq (1999).

\begin{proposition}
Let $X_{n}$ be an ARH(1) such that $E\left\Vert X_{0}\right\Vert
^{4}<+\infty,$ then the tensorized process%
\begin{equation*}
Z_{i}=X_{i}\otimes X_{i}-\Gamma
\end{equation*}
is an autoregressive process with values in $\mathcal{S}$ such that:%
\begin{equation*}
Z_{i}=R\left( Z_{i-1}\right) +u_{i}
\end{equation*}
where $R\left( S\right) =\rho S\rho^{\ast}$ and $u_{1}$ was defined at (\ref%
{AMBP_ui}). The sequence $u_{i}$ is a martingale difference with respect to
the filtration $\sigma\left( \epsilon_{i},\epsilon_{i-1},...\right) .$
\end{proposition}

\subsection{Two issues related to the general estimation problem}

\subsubsection{Identifiability}

The moment method provides the following normal equation:%
\begin{equation}
\Delta =\rho \Gamma  \label{AMBP_mom}
\end{equation}%
where 
\begin{align*}
\Gamma & =\mathbb{E}\left( X_{1}\otimes X_{1}\right) , \\
\Delta & =\mathbb{E}\left( X_{2}\otimes X_{1}\right)
\end{align*}%
are the covariance operator (resp. the cross covariance operator of order
one) of the process $\left( X_{n}\right) _{n\in \mathbb{Z}}$.

The first step consists in checking that the Yule-Walker equation (\ref%
{AMBP_mom}) correctly defines the unknown parameter $\rho.$

\begin{proposition}
\label{ident}When the inference on $\rho$ is based on the moment equation (%
\ref{AMBP_mom}), identifiability holds if $\ker\Gamma=\left\{ 0\right\} $.
\end{proposition}

The proof of this proposition is plain since taking $\widetilde{\rho }=\rho
+u\otimes v$ where $v$ belongs to the kernel of $\Gamma ,$ whenever this set
is non-empty we see that%
\begin{equation*}
\widetilde{\rho }\Gamma =\rho \Gamma +u\otimes \Gamma v=\rho \Gamma
\end{equation*}%
hence that (\ref{AMBP_mom}) holds for $\widetilde{\rho }\neq \rho $.\newline
Consequently the injectivity of $\Gamma $ is a basic assumption which can
hardly be removed and which entails that the eigenvalues of $\Gamma $ are
infinite, strictly positive. These eigenvalues will be denoted $\left(
\lambda _{i}\right) _{i\in \mathbb{N}}$ where one assumes once and for all
that the $\lambda _{i}$'s are arranged in a decreasing order with $%
\sum_{i\in \mathbb{N}}\lambda _{i}$ finite. The corresponding eigenvectors
(resp. eigenprojectors) will be denoted $\left( e_{i}\right) _{i\in \mathbb{N%
}}$ (resp. $\left( \pi _{i}\right) _{i\in \mathbb{N}}$ where $\pi
_{i}=e_{i}\otimes e_{i}$). Heuristically we should expect with (\ref%
{AMBP_mom}) at hand that $\Gamma ^{-1}$ exists to estimate $\rho $ and this
inverse will not be defined if $\Gamma $ is not one to one.

\subsubsection{The inverse problem}

Even if the identifiability is ensured estimating $\rho$ is a difficult task
due to an underlying inverse problem which stems from display (\ref{AMBP_mom}%
). The notion of inverse (or ill-posed) problem is classical in mathematical
analysis (see for instance Tikhonov, Arsenin (1977) or Groetsch (1993)). In
our framework it could be explained by claiming that equation (\ref{AMBP_mom}%
) will imply that any attempt to estimate $\rho$ will result in a highly
unstable estimate. This comes down with simple words, which will be
developed below, from the inversion of $\Gamma.$ A canonical example of an
inverse problem is the numerical inversion of an ill-conditioned matrix
(that is a matrix with eigenvalues close to zero).

The first stumbling stone comes from the fact that we cannot deduce from (%
\ref{AMBP_mom}) that $\Delta\Gamma^{-1}=\rho$. We know that a sufficient
condition for $\Gamma^{-1}$ to be defined as a linear mapping is: $%
\ker\Gamma=\left\{ 0\right\}$. Then $\Gamma^{-1}$ is an unbounded symmetric
operator on $H$. Some consequences are collected in the next proposition:

\begin{proposition}
When $\Gamma $ is injective $\Gamma ^{-1}$ may be defined. It is a linear
measurable mapping defined on a dense domain in $H,$ denoted $\mathcal{D}%
\left( \Gamma ^{-1}\right) $ and defined by:%
\begin{equation*}
\mathcal{D}\left( \Gamma ^{-1}\right) =\mathrm{Im}\Gamma =\left\{
x=\sum_{p=1}^{+\infty }x_{p}e_{p}\in H,\ \sum_{p=1}^{+\infty }\dfrac{%
x_{p}^{2}}{\lambda _{p}^{2}}<+\infty \right\} .
\end{equation*}%
This domain is dense in $H$. It is not an open set for the norm topology of $%
H$. The operator is unbounded which means that it is continuous at no point
of $\mathcal{D}\left( \Gamma ^{-1}\right) $. Besides $\Gamma ^{-1}\Gamma
=I_{H}$ but $\Gamma \Gamma ^{-1}=I_{\mathcal{D}\left( \Gamma ^{-1}\right) }$
and $\Gamma \Gamma ^{-1}$, which is not defined on the whole $H$, may be
continuously extended to $H$.
\end{proposition}

For similar reasons (\ref{AMBP_mom}) implies $\Delta\Gamma^{-1}=\rho_{|%
\mathrm{Im}\Gamma}\neq\rho$ and $\Delta\Gamma^{-1}$ may be continuously and
formally extended to the whole $H.$ In fact $\Gamma$ hence $\Gamma^{-1}$ are
unknown. However would $\Gamma$ be totally accessible we should find a way
to regularize the odd mathematical object that is $\Gamma^{-1}$. Within the
literature on inverse problems (see for instance Groetsch (1993)) one often
replaces $\Gamma^{-1}$ by a linear operator "close" to it but endowed with
additional regularity (continuity/boundedness) properties, say $%
\Gamma^{\dag} $. The Moore-Penrose pseudo inverse is an example of such an
operator but many other techniques exist. Indeed starting from 
\begin{equation*}
\Gamma^{-1}\left( x\right) =\sum_{l\in\mathbb{N}}\dfrac{1}{\lambda_{l}}%
\pi_{l}\left( x\right)
\end{equation*}
for all $x$ in $\mathcal{D}\left( \Gamma^{-1}\right) $ one may set for
instance:%
\begin{align}
\Gamma^{\dag}\left( x\right) & =\sum_{l\leq k_{n}}\dfrac{1}{\lambda_{l}}%
\pi_{l}\left( x\right)  \label{spectral-cut} \\
\Gamma^{\dag}\left( x\right) & =\sum_{l}\dfrac{1}{\lambda_{l}+\alpha_{n}}%
\pi_{l}\left( x\right)  \label{penal} \\
\Gamma^{\dag}\left( x\right) & =\sum_{l}\dfrac{\lambda_{l}}{\lambda
_{l}^{2}+\alpha_{n}}\pi_{l}\left( x\right)  \label{tikho}
\end{align}
where $k_{n}$ is an increasing and unbounded sequence of integers and $%
\alpha_{n}$ a sequence of positive real numbers decreasing to $0$. The three
operators in the display above are indexed by $n$, are all bounded with
increasing norm and are known as the spectral cut-off, penalized and
Tikhonov regularized inverses of $\Gamma$. They share the following
pointwise convergence property:%
\begin{equation*}
\Gamma^{\dag}x\rightarrow\Gamma^{-1}x
\end{equation*}
for all $x$ in $\mathcal{D}\left( \Gamma^{-1}\right) $.

In practice if $\Gamma_{n}$ is a convergent estimator of $\Gamma$ the
regularizing methods introduced below can be applied to $\Gamma_{n}$ which
is usually not invertible (see below for an example). It should be noted at
this point that the regularization for the inverse of the covariance
operator appears in the linear regression model for functional variable:%
\begin{equation*}
y=\left\langle X,\varphi\right\rangle +\epsilon
\end{equation*}
when estimating the unknown $\varphi$ (see Cardot et al. (2007)).

At last a general scheme to estimate $\rho$ may be proposed with estimates
of $\Gamma$ and $\Delta$ at hand say $\Gamma_{n}$ and $\Delta_{n}$: compute $%
\Gamma_{n}^{\dag}$ and take for the estimate and the predictor based on the
new input $X_{n+1}$ respectively:%
\begin{equation*}
\widehat{\rho}_{n}=\Delta_{n}\Gamma_{n}^{\dag}\quad\mathrm{and}\quad 
\widehat{\rho}_{n}\left( X_{n+1}\right) .
\end{equation*}
Obviously examples of such estimates are the empirical covariance and
cross-covariance operators

\begin{align*}
\Gamma _{n}& =\dfrac{1}{n}\sum_{k=1}^{n}X_{k}\otimes X_{k}, \\
\Delta _{n}& =\dfrac{1}{n-1}\sum_{k=1}^{n-1}X_{k+1}\otimes X_{k}
\end{align*}%
where the $X_{k}$'s were reconstructed by interpolation techniques.

For instance the spectral cut-off version for $\Gamma_{n}$ is 
\begin{equation*}
\Gamma_{n}^{\dag}=\sum_{l\leq k_{n}}\dfrac{1}{\widehat{\lambda}_{l}}\widehat{%
\pi}_{l}
\end{equation*}
where the eigenvalues $\widehat{\lambda}_{l}$ and the eigenprojectors $%
\widehat{\pi}_{l}$ are by-products of the functional PCA of the sample $%
X_{1},...,X_{n}.$

\begin{remark}
This inverse problem is the main serious abstract concern when infering on
the ARH model. The considerations above are moreless exposed in all the
articles dealing with it and we guess it will be of some interest to expose
and sum up this issue and some of its solutions in this monograph.
\end{remark}

\subsection{Convergence results for the autocorrelation operator and the
predictor}

As the data are of functional nature, the inference on $\rho$ cannot be
based on likelihood. Lebesgue's measure cannot be defined on
infinite-dimensional spaces. However it must be mentioned that Mourid and
Bensma\"{\i}n (2006) propose to adapt Grenander's theory of sieves
(Grenander (1981) and Geman and Hwang (1982)) to this issue. They prove
consistency in two very important cases: when $\rho$ is a kernel operator
and when $\rho$ is Hilbert-Schmidt. In the former case $\rho$ is identified
with the associated kernel $K,$ developed on a basis of trigonometric
functions along the sieve:%
\begin{equation*}
\Theta_{m}=\left\{ K\in L^{2}:K\left( t\right) =c_{0}+\sum_{k=1}^{m}\sqrt{2}%
c_{k}\cos\left( 2\pi kt\right) ,\ t\in\left[ 0,1\right] ,\
\sum_{k=1}^{m}k^{2}c_{k}^{2}\leq m\right\} .
\end{equation*}

This approach is truly original within the literature on functional data and
could certainly be extended to other problems of linear or non linear
regression.

The seminal paper dealing with the estimation of the operator $\rho$ dates
back to 1991 and is due to Bosq (1991). Several consistency results are
carried out immediately relayed by Pumo's (1992) and Mourid's (1995) PhD
thesis. Then Pumo (1998) focus on random functions with values in $C\left( %
\left[ 0,1\right] \right) $ with specific techniques. Besse and Cardot
(1996), then Besse, Cardot and Stephenson (2000) implement spline and kernel
methodology with application to climatic variations. Amongst several
interesting ideas they introduce a local covariance estimate:%
\begin{equation*}
\widehat{\Gamma}_{h_{n}}=\frac{\sum_{i=1}^{n}\left[ X_{i}\otimes X_{i}\right]
K\left( \left\Vert X_{i}-X_{n}\right\Vert /h\right) }{\sum_{i=1}^{n-1}K%
\left( \left\Vert X_{i}-X_{n}\right\Vert /h\right) }
\end{equation*}
and a local cross-covariance estimate which emphasize data close to the last
observation. This method make it possible to consider data with departures
from the stationarity assumption. This issue of the estimation of $\rho$ is
also treated in Guillas (2001) and Mas (2004).

A recent paper by Antoniadis and Sapatinas (2003) carry out wavelet
estimation and prediction in the ARH(1) model. The inverse problem is
underlined through a class of estimates stemming from the deterministic
literature on this topic. This class of estimates is compatible with wavelet
techniques and lead to consistency of the predictor. The method is applied
on the "El Nino" dataset which tends to become a benchmark for comparing the
performances of the predictions.

Ruiz-Medina et al (2007) consider the functional principal oscillation
pattern (POP) decomposition of the operator $\rho $ as an alternative to
functional PCA decomposition. They implement a Kalman filter to the
state-space equation obtained at the preceding step and derive the optimal
predictor. This original approach, illustrated by some simulations, seems to
be suited to spatial functional data as well.

Kargin and Onatski (2008) introduce the notion of predictive factor which
seems to be better suited than the PCA basis to project the data if one
really focuses on the predictor (and not on the operator itself). A double
regularization (penalization and projection) provides them with a rate of $%
O\left( n^{-1/6}\log ^{\beta }n\right) $ (where $\beta >0$) for the
prediction mean square error.

In Mas (2007) the problem of weak convergence is addressed. The main results
are given in the Theorem below:

\begin{theorem}
\label{ThMas}It is impossible for $\widehat{\rho}_{n}-\rho$ to converge in
distribution for the classical norm topology of operators. But under moment
assumptions, if $\left\Vert \Gamma^{-1/2}\rho\right\Vert _{\infty}<+\infty$
and if the spectrum of $\Gamma$ is convex then when $k_{n}=o\left( \dfrac{%
n^{1/4}}{\log n}\right) ,$%
\begin{equation*}
\sqrt{\dfrac{n}{k_{n}}}\left( \widehat{\rho}_{n}\left( X_{n+1}\right) -\rho%
\widehat{\Pi}_{k_{n}}\left( X_{n+1}\right) \right) \overset{w}{\rightarrow}%
\mathcal{G}
\end{equation*}
where $\mathcal{G}$ is a $H$-valued gaussian centered random variable with
covariance operator $\Gamma_{\epsilon}$ and $\widehat{\Pi}_{k_{n}}$ is the
projector on the $k_{n}$ first eigenvectors of $\Gamma_{n}$.
\end{theorem}

\begin{remark}
The first sentence of the Theorem above is quite surprising but is a direct
consequence of the underlying inverse problem. Finally considering the
predictor weakens the topology and has a smoothing effect on $\widehat{\rho }%
_{n}$. This phenomenon -which was exploited in Antoniadis, Sapatinas (2003)-
appears as well in the linear regression model for functional data (see
Cardot, Mas, Sarda (2007)).
\end{remark}

It should be noted that rates of convergence are difficult to obtain (see
Guillas (2001) or Kargin and Onatski (2008), Theorem 3) and rather slow with
respect to those obtained in the regression model. An exponential inequality
appears at Theorem 8.8 in Bosq (2000) but it seems that a more systematic
study of the mean square prediction error has not been carried out yet and
that optimal bounds are not available$.$

\subsubsection{Hypothesis testing}

\label{AMBP_test_ARH}

A very recent article by Horvath, Huskova and Kokoszka (2009) focuses on the
stability of the autocorrelation operator against change-point alternatives.
In fact the model (\ref{AMBP_ARH}) based on the sample $X_{1},...,X_{n}$ is
slightly modified to: 
\begin{equation*}
X_{n}=\rho_{n}\left( X_{n-1}\right) + \epsilon_{n}
\end{equation*}
and the authors test 
\begin{equation*}
H_{0}:\rho_{1}=...=\rho_{n}
\end{equation*}
against the alternative:%
\begin{equation*}
H_{A}:\text{ there exists }k^{\ast}\in\left\{ 1,...,n\right\} :\rho
_{1}=...=\rho_{k^{\ast}}\neq\rho_{k^{\ast}+1}=..=\rho_{n}.
\end{equation*}
The test is based on the projection of the process $X_{n}$ on the $p$ first
eigenvectors of the functional PCA and on an accurate approximation of the
long-run covariance matrix. The asymptotic distribution is derived by means
of empirical process techniques. The consistency of the test is obtained and
a simulation/real case study dealing with credit card transaction time
series is treated.

It turns out that Lauka\"{\i}tis and Rackauskas (2002) considered the same
sort of problem a few years sooner. They introduce a functional version of
the partial sum process of estimated residuals:%
\begin{equation*}
S\left( t\right) =\sum_{k=2}^{\left\lfloor t\right\rfloor }\left[ X_{k}-%
\widehat{\rho}\left( X_{k-1}\right) \right]
\end{equation*}
and obtain weak convergence results for its normalized version to an $H$%
-valued Wiener process. This formal theorem yields different strategies
(dyadic increment of partial sums or moving residual sums) to derive a test.

It seems however that the topic of hypothesis testing was rarely addressed
yet quite promising even if serious theoretic and technical problems appear,
once again in connection with the inverse problem mentioned earlier in this
article.

\subsection{Extension of ARH model}

Various extensions have been proposed for ARH(1) model in order to improve
the prediction performance of ARH(1) model. The first one is the natural
extension autoregressive process of order $p$ with $p > 1$, denoted ARH(p),
defined by 
\begin{equation*}
X_{n} = \rho_1 X_{n-1} + \ldots + \rho_p X_{n-p} + \epsilon_{n} .
\end{equation*}
Using the Markov representation $Y_n=\rho^{\prime}Y_{n-1} +
\epsilon_{n}^{\prime}$ where 
\begin{equation*}
\rho^{\prime}= \left[ 
\begin{array}{cccc}
\rho_1 & \rho_2 & \cdots & \rho_n \\ 
I & 0 & \cdots & 0 \\ 
\vdots &  &  & \vdots \\ 
0 &  &  & 1%
\end{array}
\right], Y_n = \left[ 
\begin{array}{c}
X_{n} \\ 
X_{n-1} \\ 
\vdots \\ 
X_{n-p+1}%
\end{array}
\right], \text{and } \epsilon^{\prime}_n = \left[ 
\begin{array}{c}
\epsilon_{n} \\ 
0 \\ 
\vdots \\ 
0%
\end{array}
\right].
\end{equation*}
and $I$ denotes the identity operator, Mourid (2003) obtain asymptotic
results of projector estimators and predictors.

Damon and Guillas (2002) introduced autoregressive Hilbertian process with
exogenous variables model, denoted ARHX(1), which intends to take into
account the dependence structure of random curves under the influence of
explanatory variables. The model is defined by the equation 
\begin{equation*}
X_{n}=\rho (X_{n-1})+a_{1}(Z_{n,1})+\ldots a_{q}(Z_{n,q})+\epsilon _{n},n\in
Z
\end{equation*}%
where $a_{1},\cdots ,a_{q}$ are bounded linear operators in $H$ and $%
Z_{n,1},\cdots ,Z_{n,q}$ are ARH(1) exogenous variables; they suppose that
the noises of the $q+1$ H$-$ valued autoregressive processes are
independent. They obtain some limit theorems, derive consistent estimators,
present a simulation study in order to illustrate the accuracy of the
estimation and compare the forecasts with other functional models.

Guillas (2002) consider a H-valued autoregressive stochastic sequence $%
(X_{n})$ with several regimes such that the underlying process $(I_{n})$ is
stationary. Under some dependence assumptions on $(I_{n})$ he proves the
existence of a unique stationary solution and state a law of large numbers
and the consistency of the covariance estimator. Following the same idea in
a recent work Mourid (2004) introduces and studies the autoregressive
process with random operators $X_{n}=\rho _{n}X_{n-1}+\epsilon _{n}$ where $%
(\rho _{n},n\in Z)$ is stationary and independent of $(\epsilon _{n})$.
Results similar to classical $ARH(1)$ are obtained.

A new model, denoted ARHD process, considering the derivative curves of an
ARH(1) model was introduced by Marion and Pumo (2004). In a recent paper Mas
and Pumo (2007) introduced and study a slightly new model: 
\begin{equation*}
X_{n}= \phi X_{n-1} + \Psi ( X_{n-1}^{\prime}) + \epsilon_{n}
\end{equation*}
where $X_n$ are random function with values in the Sobolev space $W^{2,1} =
\{u \in L^{2} [0,1], u^{\prime}\in L^{2} [0,1] \}$, $\phi$ is a compact
operator from $W$ to $W$, $\Psi$ is a compact operator from $L^{2} [0,1]$ to 
$W^{2,1}$ and $\| \phi h + \Psi h^{\prime}\| \leq \|h\|$ for $h \in W^{2,1}$%
. Convergent estimates are obtained through an original double penalization
method. Simulations on real data show that predictions are comparable to
those obtained by other classical methods based on ARH(1) modelization.
Tests on the derivative part and models with higher derivatives may be
interesting from both theoretical and practical point of view.

\subsection{Numerical aspects}

We present in this section some numerical aspects concerning the prediction
when data are curves observed at discrete points. To our knowledge the
prediction methods based on linear processes are limited to application of
ARH(1) model since tractable algorithms using general linear processes in
Hilbert spaces do not exist (see Merlev\`{e}de (1997)). However some partial
results are available for moving average processes in Hilbert spaces which
will be briefly discussed in the next section.

The literature using ARH(1) model to make prediction is various and rich and
concern different domains:

\begin{itemize}
\item Environment: Besse et al. (2000); Antoniadis and Sapatinas (2003); Mas
and Pumo (2007); Fern\'{a}ndez de Castro et al. (2005); Damon et Guillas
(2002);

\item Economy and finance: Kargin and Onatski (2008);

\item Electricity consumption: Cavallini et al. (1994);

\item Medical sciences: Marion and Pumo (2004); Glendinning and Fleet (2007)
\end{itemize}

From a technical point of view the different approaches for implementing an
ARH proceed in two steps. The first step consists in decomposing data in
some functional basis in order to reconstruct them on the whole observed
interval. Most of the methods use spline or wavelet basis and suppose that
curves belong to the Sobolev $W^{2,k}$ space of functions such that the $k$%
-th derivative is squared integrable. We invite the reader to refer to the
papers by Besse and Cardot (1996), Pumo (1998) and Antoniadis and Sapatinas
(2003) among others for detailed discussions about the use of splines and
wavelets for numerical estimation and prediction using ARH(1) model and for
the numerical results presented hereafter. The second step consists in
choosing tuning parameters required by these methods, for example the
dimension of the projection subspace for the projection estimators. A
general method used by the precedent authors is based on cross-validation
approach which gives satisfactory results in applications. Note at last that
alternatives approaches of prediction based on ARH(1) modelization are
proposed by Mokhtari and Mourid (2002) and Mourid and Bensmain (2005). In
Mokhtari and Mourid (2002) the authors use a Parzen approximation on
reproducing kernel spaces framework. Some simulation studies are presented
in the recent paper published in 2008 by the same authors.

\begin{figure}[tbp]
\centering{\includegraphics[scale=.55]{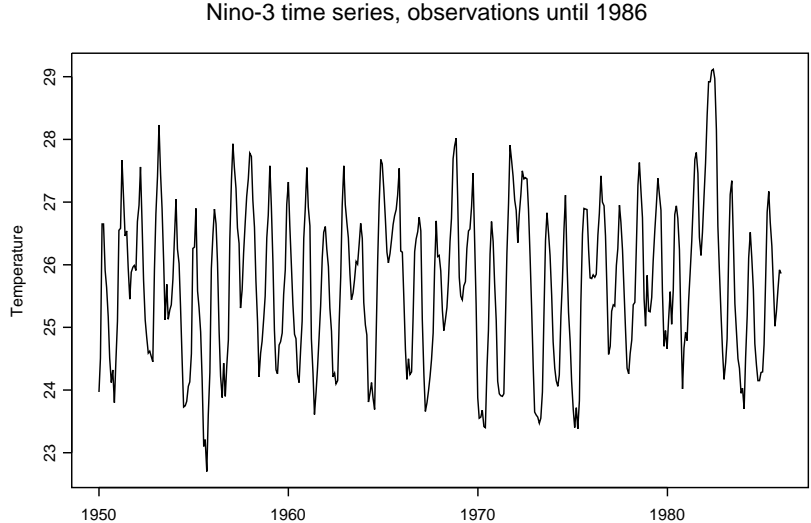}}
\caption{Monthly mean El Ni\^{n}o sea surface temperature index from January
1950 to December 1986}
\label{AMBP_sst_proc_until_1986}
\end{figure}

In order to compare methods described above we consider a climatological
time series describing the El Ni\~{n}o-Southern Oscillation (see. for
example Besse et al. (2000) or Smith et al. (1996) for a description of the
data\footnote{%
Data is freely available from
http://www.cpc.ncep.noaa.gov/data/indices/index.html}). The series gives the
monthly mean El Ni\^{n}o sea surface temperature index from January 1950 to
December 1986 and is presented in figure \ref{AMBP_sst_proc_until_1986}. We
compare the ARHD predictor with various functional prediction methods.

We compare the predictors of month temperature during 1986 knowing the data
until 1985 by two-criteria: mean-squared error (MSE) and relative
mean-absolute error (RMAE) defined by: 
\begin{equation*}
MSE=\frac{1}{12}\sum_{i=1}^{12}\left( X_{n}^{i}-\hat{X}_{n}^{i}\right)
^{2},\ RMAE=\frac{1}{12}\sum_{i=1}^{12}\frac{\left\vert X_{n}^{i}-\hat{X}%
_{n}^{i}\right\vert }{X_{n}^{i}}
\end{equation*}%
where $X_{n}^{i}$ (resp. $\hat{X}_{n}^{i}$) denotes the i$-$th month
observation (resp. prediction). The two-criteria for various functional
predictors are given in Table \ref{AMBP_sst_pred_1986}. 
Results show that the best method are Wavelet II (one of the wavelet
approaches proposed in Antoniadis and Sapatinas (2003)) and spline smoothing
ARH(1). Globally the predictors obtained using ARH(1) model are better and
numerically faster than the classical SARIMA $(0,1,1)\times (1,0,1)_{12}$
model (the best SARIMA model based on classical criteria).

\begin{table}[ptb]
\begin{center}
\begin{tabular}{|p{7cm}|c|c|}
\hline\hline
Prediction method & MSE & RMAE (\%) \\ \hline\hline
Wavelet II & 0.063 & 0.89 \\ 
Splines & 0.065 & 0.89 \\ 
ARHD & 0.167 & 1.25 \\ 
ARH(1): linear spline & 0.278 & 2.4 \\ 
SARIMA & 1.457 & 3.72 \\ \hline\hline
\end{tabular}%
\end{center}
\caption{Mean Squared Error (MSE) and RMAE errors for prediction of El Ni\^{n%
}o index during 1986}
\label{AMBP_sst_pred_1986}
\end{table}


\section{Perspectives}

In the precedent sections we insisted on two important statistical problems
concerning H linear processes. The first discussed in \S \ref%
{AMBP_LHP_perspectives} and in relation with inference on general or
generalized linear processes. The estimation with the aim to make
predictions with such processes seems to arise difficult technical problems.
Some new results in this direction are obtained recently by Bosq (2006) by
introducing the moving average process of order $q\geq 1$, $MAH(q)$. Some
partial consistency results for the particular process MAH(1) are presented
in a paper by Turbillon et al. (2008). A MAH(1) is a $H$ valued process
satisfying the equation $X_{t}=\epsilon _{t}+\ell (\epsilon _{t-1})$ where $%
\ell $ is a compact operator and $(\epsilon _{t})$) a strong white noise. It
is simple from (\ref{AMBP_cond_inv_lph}) to show that this process is
invertible when the condition $\Vert \ell \Vert <1$. The difficulty in
estimating $\ell $ as for the real valued MA processes stems from the fact
that the moment equation is not linear conversely to the ARH(1) process.
Under mild conditions Turbillon et al. (2008) propose two types of
estimators for $\ell $ and give consistency results.

The second direction concerns the ARH(1) model and his extensions. In \S \ref%
{AMBP_test_ARH} we recall some serious theoretical and technical problems
with the topic of hypothesis testing. But the problem is very important in
particular from a practical point of view. As an example let us consider the
ARHD model and the test addressing the significance of the derivative in the
model. Another issue may be the characterization of real valued processes
allowing an ARH representation or more generally linear processes. While
some examples exists admitting an ARH or MAH representation (see the book by
Bosq (2000)) a general approach to recognize real processes allowing such a
representation is an issue for future works.

The above questions are important from a theoretical point of view, in
particular for the research in statistics. For the people who analyze data
that are discretized curves it's more and more necessary to dispose of
analogue description tools as for the ARMA(p,q) real valued processes. In
this direction a work by Hyndman and Shang (2008) for visualizing functional
data and identifying functional outliers is an example.

\bigskip \textbf{Acknowledgement}. The authors thank Fr\'{e}d\'{e}ric
Ferraty, Yves Romain and the whole group STAPH for initiating this work as well as for
permanent and fruitful collaboration and are grateful to Professor Denis
Bosq for helpful discussions and pointing out recent articles about
functional linear processes.


\section*{References}

\noindent\textsc{Adams R.A., Fournier J.J.F.} (2003). \textit{Sobolev spaces}%
. Academic Press, 2nd ed.

\noindent\textsc{Antoniadis A., Sapatinas T.} (2003). Wavelet methods for
continuous-time prediction using representations of autoregressive processes
in Hilbert spaces. \textit{J. of Multivariate Analysis}, \textbf{87} 133-158.

\noindent\textsc{Besse, P. et Cardot, H.} (1996). Approximation spline de la
pr\'{e}vision d'un processus fonctionnel autor\'{e}gressif d'ordre 1. 
\textit{Canad. J. Statist}, \textbf{24} 467-487.

\noindent\textsc{Besse, P., Cardot, H. et Stephenson, D.} (2000).
Autoregressive forecasting of some climatic variations. \textit{Scand. J.
Statist}, \textbf{27} 673-687.

\noindent\textsc{Bosq, D.} (1991). Modelization, nonparametric estimation
and prediction for continuous time processes. in: \textit{Roussas (Ed), Nato
Asi Series C}, \textbf{335} 509-529.

\noindent\textsc{Bosq, D.} (1996). Limit theorems for banch-valued
autoregressive processes. Applications to real continuous time processes. 
\textit{Bull. Belg. Math. Sc}, \textbf{3} 537-555.

\noindent\textsc{Bosq, D.} (1999). Repr\'{e}sentation autor\'{e}gressive de
l'op\'{e}rateur de covariance empirique d'un ARH(1). Applications. \textit{%
C.R. Acad.Sci.}, \textbf{329} S\'{e}r. I 531-534.

\noindent\textsc{Bosq, D.} (2000). \textit{Linear processes in function
spaces}. Lectures notes in statistics, Springer Verlag.

\noindent\textsc{Bosq, D.} (2007). General linear processes in Hilbert
spaces and prediction. \textit{J. Stat. Planning and Inference}, \textbf{137}
879-894.

\noindent \textsc{Bosq, D., Blanke D.} (2007). \textit{Inference and
prediction in large dimensions}. Wiley series in probability and statistics,
John Wiley and Sons, Dunod.

\noindent \textsc{Bosq, D.} (2009). Tensor products of functional ARMA
processes, submitted article.

\noindent\textsc{Brockwell P. and Davis A.} (1991). \textit{Time series:
Theory and methods}. Springer Verlag.

\noindent\textsc{Cardot H, Mas A., Sarda P.} (2007) CLT in functional linear
regression models. Probability Theory and Related Fields, \textbf{138}
325-361.

\noindent\textsc{Cavallini, A., Montanari G.C., Loggini, M., Lessi, O.,
Cacciari M.,} (1994). Nonparmetric prediction of harmonic levels in
electrical networks. in: \textit{Proceedings of IEEE ICHPS VI. Bologna},
165-171.

\noindent\textsc{Damon, J., Guillas, S.} (2002). The inclusion of exogenous
variables in functional autoregressive ozone forecasting. \textit{%
Environmetrics}, \textbf{13} (7) 759-774.

\noindent\textsc{Damon, J., Guillas, S.} (2005). Estimation and Simulation
of Autoregressive Hilbertian Processes with Exogenous Variables. \textit{%
Stat. Infer. for Stoch. Proc.}, \textbf{8} (2) 185-204.

\noindent\textsc{Dedecker, J., Merlevede, F.} (2003). The conditional
central limit theorem in Hilbert spaces. \textit{Stochastic Process. Appl.}, 
\textbf{108} 229-262.

\noindent\textsc{Dedecker, J., Merlevede, F.} (2007). Convergence rates in
the law of large numbers for Banach valued dependent variables. \textit{%
Teor. Veroyatnost. i Primenen} \textbf{52} 562-587.

\noindent\textsc{Dembo A., Zeitouni O.} (1993). \textit{Large deviations
techniques and applications}. Jones and Bartlett, London.

\noindent\textsc{Dunford, N., Schwartz, J.T.} (1988). \textit{Linear
Operators}, Vol. I \& II. Wiley Classics Library.

\noindent\textsc{Fern\'{a}ndez de Castro, B., Guillas, S., Gonz\'{a}lez
Manteiga, W.} (2005). Functional Samples and Bootstrap for Predicting Sulfur
Dioxide Levels. \textit{Technometrics} \textbf{47} (2) 212-222.

\noindent\textsc{Ferraty F., Vieu P.} (2006). \textit{Nonparametric
functional data analysis. Theory and Practice}. Springer-Verlag, New-York.

\noindent\textsc{Geman, S., Hwang, C.R.} (1982). Nonparametric maximum
likelihood estimation by the method of sieves. \textit{Ann. Statist.} 
\textbf{10} (2) 401-414.

\noindent\textsc{Glendinning, R.H., Fleet, S.L.} (2007). Classifying
functional time series. \textit{Signal Processing}, \textbf{87} (1) 79-100.

\noindent\textsc{Gohberg, I., Goldberg, S., Kaashoek,M.A.} (1991). \textit{%
Classes of linear operators Vol I \& II. Operator Theory: advances and
applications}. Birkha\"{u}ser Verlag.

\noindent\textsc{Grenander, U.} (1981). \textit{Abstract Inference}. Wiley,
New York.

\noindent\textsc{Groetsch, C.} (1993). \textit{Inverse Problems in the
Mathematical Sciences}. Vieweg, Wiesbaden.

\noindent\textsc{Guillas, S.} (2001). Rates of convergence of
autocorrelation estimates for autoregressive Hilbertian processes. \textit{%
Statist. Probab. Lett.}, \textbf{55} (3) 281-291.

\noindent\textsc{Guillas, S.} (2002) Doubly stochastic Hilbertian processes. 
\textit{J. Appl. Probab.}, \textbf{39} (3) 566-580.

\noindent\textsc{Horvat L., Huskova M, Kokoszka P.} (2009). Testing the
stability of the functional autoregressive process. \textit{J. of
Multivariate Analysis}, \textbf{to appear}.

\noindent\textsc{Hyndman R., Shang H.L.} (2008). Bagplots, boxplots and
outlier detection for Functional Data, in: \textit{Functional and
Operatorial Statistics}, 201-208, Physica-Verlag Heidelberg.

\noindent\textsc{Kargin,V., Onatski A.} (2008). Curve forecasting by
functional autoregression. \textit{J. of Multivariate Analysis}, \textbf{99}
2508-2526.

\noindent\textsc{Labbas, A., Mourid, T.} (2002). Estimation and prediction
of a Banach valued autoregressive process. \textit{C. R. Acad. Sci. Paris}, 
\textbf{335} (9) 767-772.

\noindent\textsc{Laukaitis A., Rackauskas A.} (2002). Functional data
analysis of payment systems. \textit{Nonlinear Analysis: Modeling and Control%
}, \textbf{7} 53-68.

\noindent\textsc{Ledoux, M., Talagrand M.} (1991). \textit{Probability in
Banach spaces - Isoperimetry and Processes}, Springer Verlag.

\noindent\textsc{Maltsev V.V., Ostrovskii E.I.} (1982): Central limit
theorem for stationary processes in Hilbert space. \textit{Theor. Prob. and
its Applications}, \textbf{27} 357-359.

\noindent\textsc{Marion J.M., Pumo B.} (2004). Comparaison des mod\`{e}les
ARH(1) et ARHD(1) sur des donn\'{e}es physiologiques, \textit{Annales de
l'ISUP}, \textbf{48} (3) 29-38.

\noindent\textsc{Mas, A.} (2002). Weak convergence for the covariance
operators of a Hilbertian linear process. \textit{Stoch Process. App}, 
\textbf{99} 117-135.

\noindent\textsc{Mas, A.} (2004). Consistance du pr\'{e}dicteur dans le mod%
\`{e}le ARH(1): le cas compact. \textit{Annales de l'Isup}, \textbf{48}
39-48.

\noindent\textsc{Mas, A.} (2007). Weak convergence in the functional
autoregressive model. \textit{J. of Multivariate Analysis}, \textbf{98}
1231-126.

\noindent\textsc{Mas A., Menneteau L.} (2003a). Perturbation appraoch
applied to the asymptotic study of random operators. \textit{Progress in
Probability}, \textbf{55} 127-134.

\noindent\textsc{Mas A., Menneteau L.} (2003b). Large and moderate
deviations for infinite-dimensional autoregressive processes. \textit{J. of
Multivariate Analysis,} \textbf{87} 241-260.

\noindent\textsc{Mas, A., Pumo, B.} (2007). The ARHD process. \textit{J. of
Statistical Planning and Inference}, \textbf{137} 538-553.

\noindent\textsc{Menneteau, L} (2005) Some laws of the iterated logarithm in
Hilbertian autoregressive models. \textit{J. of Multivariate Analysis}, 
\textbf{92} 405-425.

\noindent\textsc{Merlev\`{e}de, F.} (1995). Sur l'inversibilit\'{e} des
processus lin\'{e}aires \`{a} valeurs dans un espace de Hilbert. \textit{C.
R. Acad. Sci. Paris}, \textbf{321} S\'{e}rie I, 477-480.

\noindent\textsc{Merlev\`{e}de, F.} (1996). Central limit theorem for linear
processes with values in Hilbert space. \textit{Stoch. Proc. and their
Applications}, \textbf{65} 103-114.

\noindent\textsc{Merlev\`{e}de, F.} (1997). R\'{e}sultats de convergence
presque s\^{u}re pour l'estimation et la pr\'{e}vision des processus lin\'{e}%
aires hilbertiens. \textit{C. R. Acad. Sci. Paris}, \textbf{324} S\'{e}rie I
573-576.

\noindent\textsc{Merlev\`{e}de, F., Peligrad M., Utek, S.} (1997). Sharp
conditions for the CLT of linear processes in Hilbert space. \textit{J. of
Theoritical Probability}, \textbf{10} (3) 681-693.

\noindent\textsc{Mokhtari, F., Mourid, T.} (2002) Prediction of
autoregressive processes via the reproducing kernel spaces. \textit{C. R.
Acad. Sci. Paris, Ser.}, \textbf{334} 65-70.

\noindent\textsc{Mourid T.} (1995). Contribution \`{a} la statistique des
processus autor\'{e}gressifs \`{a} temps continu. PHD. Thesis, Univ. Paris
VI.

\noindent\textsc{Mourid, T.} (2002). Estimation and prediction of functional
autoregressive processes. \textit{Statistics}, \textbf{36} (2) 125-138.

\noindent\textsc{Mourid, T.} (2004). The hilbertian autoregressive process
with random oparator (in french). \textit{Annales de l'ISUP}, \textbf{48}
(3) 79-86.

\noindent\textsc{Mourid, T., Bensmain, N.} (2006). Sieves estimator of the
operator of a functional autoregressive process. \textit{Stat. and Probab.
Letters}, \textbf{76} (1) 93-108.

\noindent\textsc{Nazarova A.N.} (2000). Normal approximation for linear
stochastic processes and random fields in Hilbert space. Math. Notes, 68
363-369.

\noindent\textsc{Pumo B.} (1992). Estimation et pr\'{e}vision de processus
autor\'{e}gressifs fonctionnels. Applications aux processus \`{a} temps
continu. \textit{PHD Thesis}, Univ.Paris VI.

\noindent\textsc{Pumo B.} (1998). Prediction of continuous time processes by 
$C[0,1]$-valued autoregressive process. \textit{Statist. Infer. for Stoch.
Processes}, \textbf{3} (1) 297-309.

\noindent\textsc{Ramsay J.O.}(2000). Differential equation models for
statistical functions. \textit{Canad. J. Statist}, \textbf{28} 225-240.

\noindent\textsc{Ramsay J.O., Silverman B.W.} (1997). \textit{Functional
Data Analysis}, Springer Series in Statistics.

\noindent \textsc{Ramsay J.O., Silverman B.W.} (2002). \textit{Applied
Functional Data Analysis}. Springer Series in Statistics, Springer-Verlag,
New-York.

\noindent \textsc{Ruiz-Medina M.D., Salmeron R., Angulo J.M.} (2007). Kalman
filtering from POP-based diagonalization of ARH(1), \textit{Comp. Statist.
Data Anal.,} \textbf{51}, 4994--5008.

\noindent\textsc{Smith T.M., Reynolds R.W., Livezey R.E., Stokes D.C.}
(1996). Reconstruction of Historical Sea Surface Temperatures Using
Empirical Orthogonal Functions. \textit{Journal of Climate}, \textbf{9} (6)
1403-1420.

\noindent\textsc{Tikhonov A.N., Arsenin V.Y.} (1977). \textit{Solutions of
ill-posed problems}. V.H. Winstons and sons, Washington.

\noindent\textsc{Turbillon, S., Bosq, D., Marion, J.M., Pumo, B.} Parameter
estimation of moving averages in Hilbert spaces, \textit{C. R. Acad. Sci.
Paris, Ser.}, \textbf{346} 347-350.

\noindent\textsc{Wong, E.} (1966). Some results concerning the
zero-crossings of Gaussian noise, \textit{SIAM J. Appl. Math.}, \textbf{14}
6 1246-1254.

\noindent\textsc{Ziemer W.P.} (1989). \textit{Weakly differentiable
functions. Sobolev spaces and functions of bounded variations}. Graduate
Text in Mathematics 120, Springer-Verlag, New-York.

\end{document}